\documentclass[11pt]{article}
\usepackage{amssymb}
\usepackage{graphicx}
\usepackage{color}
\usepackage{amsmath}

\newcommand{\R}{\Bbb R}

\def\of{\, \mbox{\small{$\circ$}} \:}

\newcommand{\SO}{\mbox{\rm SO}}

\newcommand{\inv}{^{-1}}

\newcounter{investigation}
\newenvironment{rem}[2]{ \medskip
\begin{center}\begin{minipage}{5in}\refstepcounter{investigation}
\label{#2}
\noindent {\bf #1 (\Alph{investigation}).}}
{\end{minipage}\end{center}
\medskip}

\begin{document}
\bibliographystyle{plain}
\begin{center}
  {\Large An Elementary Introduction to the Hopf Fibration}\\

\end{center}

\begin{tabular}{ll}
  Author: & David W. Lyons\\
  Original publication:& April 2003\\
  Additional commentary:& 3 December 2022
\end{tabular}

The article in the following pages was published April 2003 in Volume
76, Number 2, of {\em Mathematics Magazine}, pp.87--98. This preprint is
slightly different from the published version, but only in a few minor
details.

\subsection*{Errata in the published version}

A sample computation in quaternion multiplication appears on p.89 of the
published version. The middle line on the right hand side contains two
typographical errors. A {\em correct} version of the same computation
appears on p.4 of this preprint, where the erroneous
middle line is omitted. Here is the correction for the published version.
\begin{align*}
  (3+2j)(1-4i+k) &= 3 -12i +3k +2j -8ji +2jk \hspace*{.2in}\mbox{(distributing)}\\
  &= 3-12i + 3k+2j+8k+2i \hspace*{.2in}\mbox{(applying relations)}\\
  &= 3 - 10i  + 2j +11k \hspace*{.2in}\mbox{(combining terms)}
\end{align*}

A second error appears in Figure 10, on p.96 of the published version,
where the label ``$S^3$'' for the 3-sphere appears erroneously with a
stray minus sign, as ``$S^{-3}$''. This error does not appear in this
preprint.

\subsection*{Hopf's original map}

Equation~(1) (p.87 in the published version, p.2 in this preprint)
differs from Hopf's original map in his 1931
paper by permutations of the input and output variables. For the sake of
completeness, here is Hopf's original map.
$$(a,b,c,d) \to \left(2(ac+bd),2(bc-ad),a^2+b^2-c^2-d^2\right)$$

\thispagestyle{empty}

\newpage

\setcounter{page}{1}
\thispagestyle{empty}

\begin{center}
  {\Large An Elementary Introduction to the Hopf Fibration}\\
  {published: April 2003}
\end{center}

\begin{flushright}
David W.~Lyons\\
Department of Mathematical Sciences \\
Lebanon Valley College\\
101 N. College Avenue\\
Annville, PA 17003\\
email: lyons@lvc.edu
\end{flushright}

\vskip .25 in

\subsection*{Introduction}

The Hopf fibration, named after Heinz Hopf who studied it in a 1931
paper~\cite{hopforig}, is an important object in mathematics and
physics.  It was a landmark discovery in topology and is a fundamental
object in the theory of Lie groups.  The Hopf fibration has a wide
variety of physical applications including magnetic
monopoles~\cite{nakahara}, rigid body mechanics~\cite{marsden} and
quantum information theory~\cite{mosseridandoloff}.  


Unfortunately, the Hopf fibration is little known in the undergraduate
curriculum, in part because presentations usually assume background in
abstract algebra or manifolds.  However, this is not a necessary
restriction.  We present in this article an introduction to the Hopf
fibration that requires only linear algebra and analytic geometry.  In
particular, no vector calculus, abstract algebra or topology is needed.
Our approach uses the algebra of quaternions and illustrates some of the
algebraic and geometric properties of the Hopf fibration.  We explain
the intimate connection of the Hopf fibration with rotations of 3-space
that is the basis for its natural applications to physics.

We deliberately leave some of the development as exercises, called
``Investigations,'' for the reader.  The Investigations contain key
ideas and are meant to be fun to think about.  The reader may also take
them as statements of facts that we wish to assume without
interrupting the narrative.

\subsection*{Hopf's mapping}

The {\em standard unit $n$-sphere} $S^n$
is the set of points $(x_0,x_1,\ldots,x_{n})$ in $\R^{n+1}$ that
satisfy the equation
$$x_0^2+x_1^2+\cdots+x_{n}^2=1.$$ Geometrically stated, $S^n$ is the set
of points in $\R^{n+1}$ whose distance from the origin is 1.  Thus the
1-sphere $S^1$ is the familiar unit circle in the plane, and the
2-sphere $S^2$ is the surface of the solid unit ball in
3-space.
The thoughtful reader may wonder what higher dimensional spheres look
like.  We address this issue at the end of this article, where we
explain how {\em stereographic projection} is used to ``see'' inside 
$S^3$.

The Hopf fibration is the mapping
$h\colon S^3\to S^2$ defined by
\begin{equation}\label{hopfdef}
h(a,b,c,d)=(a^2+b^2-c^2-d^2,2(ad+bc),2(bd-ac)).
\end{equation}
To be historically precise, Hopf's original formula differs from that
given here by a reordering of coordinates.  We use this altered version
to be consistent with the quaternion approach explained later in this
article.  It is easy to check that the squares of the three coordinates
on the right hand side sum to $(a^2+b^2+c^2+d^2)^2=1$, so that the image
of $h$ is indeed contained in $S^2$.

What problem was Hopf trying to solve when he invented this map?  And
how can one see any connection with physical rotations, as we have
claimed?

The work in Hopf's paper~\cite{hopforig} was an early achievement in the
modern subject of {\em homotopy theory}.  In loose terms, homotopy seeks
to understand those properties of a space that are not altered by
continuous deformations.  One way to detect these properties 
in an unknown space $X$ is to compare $X$ with a well understood space $Y$ via
the set of all continuous maps $Y\to X$.  Two maps whose images can be
continuously deformed from one to the other are considered equivalent.
Knowing something about $Y$ and also about the set of homotopically
equivalent maps from $Y$ to $X$ helps us understand $X$.  This seemingly
indirect method provides a powerful way to analyze spaces.

Ironically, one of the most intractable problems in homotopy theory is
to determine the homotopy equivalence classes of maps $Y\to X$ when $X$
and $Y$ are both spheres and the dimension of $X$ is smaller than the
dimension of $Y$.  Many individual cases for particular pairs of
dimensions of $X$ and $Y$ are understood, but there remain interesting
unsolved problems. Hopf's map $h\colon S^3 \rightarrow S^2$ was a
spectacular breakthrough in this area.  We cannot give the full story of
this discovery here, but we can explain the Hopf fibration in a
geometric way that indicates its connection to rotations.

\subsection*{Rotations and quaternions}

\begin{figure}[h]
\begin{center}
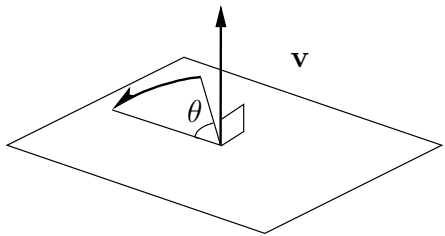
\end{center}
\vspace*{-.2in}\caption{\small A rotation in $\R^3$ is specified by an
angle $\theta$
and a vector ${\mathbf v}$ giving the axis.}
\label{vectangrot}
\end{figure}

First, notice that a rotation about the origin in $\R^3$ can be
specified by giving a vector for the axis of rotation and an angle of
rotation about the axis.  We make the convention that the rotation will be
counterclockwise for positive angles, and clockwise for negative angles,
when viewed from the tip of the vector (see Figure~\ref{vectangrot}).

The specification of a rotation by an axis vector and an angle is far
from unique.  The rotation determined by the vector ${\mathbf v}$ and
the angle $\theta$ is the same as the rotation determined by the pair
$(k{\mathbf v},\theta + 2n\pi)$, where $k$ is any positive scalar and
$n$ is any integer.  The pair $(-{\mathbf v},-\theta)$ also determines
the same rotation.  Nonetheless, we see that four real numbers are
sufficient to specify a rotation: three coordinates for a vector and one
real number to give the angle.  This is far fewer than the {\em nine}
entries of a $3\times 3$ orthogonal matrix we learn to use in linear
algebra.  In fact, we can cut the number of parameters needed to
specify a rotation from four to three, for example, by giving an axis
vector whose length determines the angle of rotation.  However, we shall
not pursue that here; it is the 4-tuple approach that turns out to be
practical.  Is there an efficient way to work with 4-tuples of real
numbers to do practical calculations with rotations?  Here is a sample
basic question.

\begin{rem}{Investigation}{composerotations}
Given geometric data (axes and angles) for two rotations, how do you
determine the axis and angle for their composition?  (The composition of
two rotations is the motion obtained by performing first one rotation,
then the other.  Order counts!)  Suggestion: think about this
investigation long enough to realize that it is difficult, or at least
tedious, if you restrict yourself to matrix methods; then revisit after
you do Investigation~(E) below.
\end{rem}

The problem of finding a convenient algebraic method for computing with
rotations led William Rowan Hamilton to invent the {\em quaternions} in
the mid-19th century.  The discovery of quaternions, and Hamilton's life
in general, is a fascinating bit of history.  For further reading,
see references~\cite{hankins} and~\cite{odonnell}.  For an exposition of the
rotation problem in Investigation~(A)
and its solution, beyond what appears in this section,
see~\cite{kuipers}, \S~6.2 ff.

Hamilton was inspired by the solution to the analogous problem in two
dimensions: rotations of the plane about the origin can be encoded by
unit length {\em complex numbers}.  The angle of a rotation is the same
as the angle made by its corresponding complex number, thought of as a
vector in $\R^2$, with the positive real axis.  The composition of
rotations corresponds to the multiplication of the corresponding complex
numbers.  Hamilton tried for years to make an algebra of rotations in
$\R^3$ using ordered triples of real numbers.
One day he realized he could achieve his goal using 4-tuples.  Here is
his invention.

As a set (and as a vector space) the set of quaternions is identical to
$\R^4$.  The three distinguished coordinate vectors $(0,1,0,0)$,
$(0,0,1,0)$ and $(0,0,0,1)$ are given the names $i$, $j$ and $k$,
respectively.  The vector $(a,b,c,d)$ is written $a+bi+cj+dk$ when
thought of as a quaternion.  The number $a$ is referred to as the {\em
real part} and $b$, $c$ and $d$ are called the {\em $i$, $j$ and $k$
parts}, respectively.  Like real and complex numbers, quaternions can be
{\em multiplied}.  The multiplication rules are encapsulated by the
following relations.
$$i^2=j^2=k^2=-1$$
$$ij=k \hspace{.25in} jk=i \hspace{.25in} ki=j$$
The elements $i$, $j$ and $k$ do {\em not} commute.  Reversing the
left-right order changes the sign of the product.
$$ji=-k \hspace{.25in} kj=-i \hspace{.25in} ik=-j$$
Here is a sample multiplication.
\begin{eqnarray*}
(3+2j)(1-4i+k) &=& 3 -12i +3k +2j -8ji +2jk 
\hspace{.2in}\mbox{(distributing)}\\
&=& 3 - 10i  + 2j +11k
\hspace{.2in}\mbox{(applying relations and combining terms)}
\end{eqnarray*}
Similar to the complex numbers, the {\em conjugate} of a quaternion
$r=a+bi+cj+dk$, denoted $\overline{r}$, is defined to be
$\overline{r}=a-bi-cj-dk$.  The {\em length} or {\em norm} of a
quaternion $r$, denoted $||r||$, is its length as a vector in $\R^4$.
(The term {\em norm}, when applied to quaternions, is also commonly used
to refer to the {\em square} of the Euclidean norm defined here.)  The
formula for the norm of $r=a+bi+cj+dk$ is
$||r||=\sqrt{a^2+b^2+c^2+d^2}$.

\begin{rem}{Investigation}{quaternions}
What algebraic properties do the quaternions share with the real or
complex numbers?  How are they different?  Show that quaternion
multiplication is associative but noncommutative.  Associativity means
that
$$p(qr)=(pq)r$$ for all quaternions $p$, $q$ and $r$.  Another formula
for the norm of $r=a+bi+cj+dk$ is $||r||=\sqrt{r\overline{r}}$.  Norm
has the property $||rs||=||r||\;||s||$ for all quaternions $r$ and $s$.
Because of this, multiplying two unit length quaternions yields another
unit length quaternion.  The set of unit length quaternions, viewed as
points in $\R^4$, is the 3-sphere $S^3$.  Each nonzero quaternion $r$ has
a {\em multiplicative inverse}, denoted $r^{-1}$, given by
$$r^{-1}=\frac{\overline{r}}{||r||^2}.$$  
When $r$ is a unit quaternion, $r\inv$ is the same as $\overline{r}$.
(For these and other details about quaternion algebra,
see~\cite{kuipers} Ch.~5.)
\end{rem}

Here is how a quaternion $r$ determines a linear mapping
$R_r\colon \R^3\to\R^3$.  To a point $p=(x,y,z)$ in 3-space, we
associate a quaternion $xi+yj+zk$ which we will also call $p$
(a quaternion whose real part is zero is called {\em pure}).  The
quaternion product $rpr\inv$ is also pure, that is, of the form
$x'i+y'j+z'k$, and hence can be thought of as a point $(x',y',z')$ in
3-space.  We define the mapping $R_r$ by
\begin{equation}\label{quatrotdef}
R_r(x,y,z)=(x',y',z').
\end{equation}

\begin{figure}[h]
\begin{center}
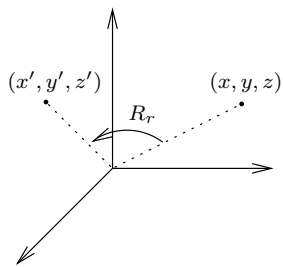
\end{center}
\vspace*{-.2in}\caption{\small A nonzero quaternion $r$ gives rise to a rotation $R_r$ in $\R^3$.}
\end{figure}

\begin{rem}
{Investigation}{quatmapislinear} Is the mapping $R_r$ described in the
previous paragraph indeed a linear map?  Verify that this is the case.
Moreover, show that the map determined by any nonzero real scalar
multiple of $r$ is equal to $R_r$, that is, show that $R_{kr} = R_r$ for
any quaternion $r$ and any nonzero real number $k$.  Show that when
$r\neq 0$, $R_r$ is invertible with inverse $(R_r)\inv = R_{(r\inv)}$.
\end{rem}

From the ``moreover'' statement in the previous Investigation, whenever
$r\neq 0$, we are free to choose $r$ to have norm 1 when working with the
map $R_r$, and we shall do so since this makes the analysis simpler.

For $r\neq 0$, it turns out that $R_r$ is a rotation of $\R^3$.  The
axis and angle of the rotation $R_r$ are elegantly encoded in the four
coordinates $(a,b,c,d)$ in the following way, when $r$ is a unit
quaternion.  If $r=\pm 1$, it is easy to see that $R_r$ is the identity
mapping on $\R^3$.  Otherwise, $R_r$ is a rotation about the axis
determined by the vector $(b,c,d)$, with angle of rotation
$\theta=2\cos^{-1}(a)=2\sin^{-1}(\sqrt{b^2+c^2+d^2})$.  To appreciate
how nice this is, have a friend write down a $3\times 3$ orthogonal
matrix, say, with no zero entries; now quickly find the axis and angle
of rotation!

The facts stated in the preceding paragraph are not supposed to be
obvious.  Here is a sequence of exercises that outline the proof.  For
a detailed discussion see~\cite{kuipers} \S~5.15.

\begin{rem}{Investigation}{quatrotation}
How does a unit quaternion encode geometric information about its
corresponding rotation?  Let $r=a + bi + cj + dk$ be a unit quaternion.
Verify that if $r=\pm 1$, then $R_r$ defined above is the identity
mapping.  Otherwise, show that $R_r$ is the rotation about the axis
vector $(b,c,d)$ by the angle
$\theta=2\cos^{-1}(a)=2\sin^{-1}(\sqrt{b^2+c^2+d^2})$, as follows.
\begin{enumerate}
\item Show that $R_r$ preserves norm, i.e., that $||R_r(p)|| = ||p||$ for
any pure quaternion $p = xi + yj + zk$.  (This follows from the fact
that the norm of a quaternion product equals the product of the norms.)
\item Show that the linear map $R_r$ has eigenvector $(b,c,d)$ with eigenvalue 1.
\item Here is a strategy to compute the angle of rotation.  Choose a
vector ${\mathbf w}$ perpendicular to the eigenvector $(b,c,d)$.  This
can be broken down into two cases: if at least one of $b,c$ is nonzero,
we may use ${\mathbf w} = ci-bj$.  If $b=c=0$, we may use ${\mathbf
w}=i$.  Now compute the angle of rotation by finding the angle between
the vectors ${\mathbf w}$ and $R_r{\mathbf w}$ using the following
formula from analytic geometry, where the multiplication in the
numerator on the right hand side is the dot product in $\R^3$.
$$\cos \theta = \frac{{\mathbf w}\cdot R_r{\mathbf w}}{||{\mathbf
w}||^2}$$ In all cases the right hand side equals
$a^2-b^2-c^2-d^2=2a^2-1$.  Now apply a half-angle identity to get
$a=\cos(\theta/2)$.
\end{enumerate}
\end{rem}

Here is the fact that illustrates how Hamilton accomplished his goal to
make an algebra of rotations.

\begin{rem}
{Investigation}{composequat}
Let $r$ and $s$ be unit quaternions.  Verify that
$$R_r \of R_s = R_{rs}.$$
In words, the {\em composition} of rotations can be accomplished by the
{\em multiplication} of quaternions.  Now go back and try
Investigation~(A).
\end{rem}

The next Investigation is appropriate for a student who has some
experience with groups, or could be a motivating problem for an
independent study in the basics of group theory. (For an excellent
introduction to group theory with a geometric point of view,
see~\cite{armstrong}.)

\begin{rem}{Investigation}{sthreeisgrp}
The set $S^3$ with the operation of quaternion multiplication satisfies
the axioms of a {\em group}.  The set of rotations in 3-space, with the
operation of composition, is also a group, called $\SO(3)$.  The map
$\varphi\colon S^3 \rightarrow \SO(3)$ given by $r\mapsto R_r$ is a {\em group
homomorphism}.  Each rotation $R$ in $\SO(3)$ can be written in the
form $R=R_r$ for some $r\in S^3$ (i.e., the map $\varphi$ is {\em surjective}),
and each rotation $R_r$ has precisely two preimages in $S^3$, namely
$r$ and $-r$.  The {\em kernel} of $\varphi$ is the subgroup $\{1,-1\}$,
and we have an {\em isomorphism} of groups
$$S^3/\{1,-1\} \approx SO(3).$$
\end{rem}

\subsection*{The 3-sphere, Rotations, and the Hopf Fibration}

We now give a reformulation of the Hopf map in terms of quaternions.
First, fix a distinguished point, say, $P_0=(1,0,0)$, on $S^2$.
Given a point $(a,b,c,d)$ on $S^3$, let $r=a+bi+cj+dk$ be the
corresponding unit quaternion.  The quaternion $r$ then defines a rotation $R_r$ of
3-space given by~(\ref{quatrotdef}) above.  Then the Hopf
fibration is defined by
\begin{equation}\label{hopfdef2}
r\mapsto R_r(P_0) = ri\overline{r}.
\end{equation}
\begin{rem}{Investigation}{hopffibviarot}
Verify that the two formulas~(\ref{hopfdef}) and~(\ref{hopfdef2}) for
the Hopf fibration are equivalent.
\end{rem}

\begin{figure}[h]
\begin{center}
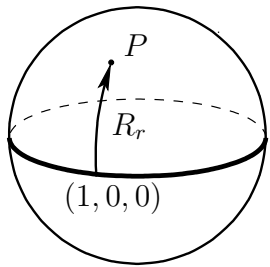
\end{center}
\vspace*{-.2in}\caption{\small The unit quaternion $r$ moves $(1,0,0)$ to $P$ via $R_r$.
The Hopf map takes $r$ to $P$.}
\end{figure}


Consider the point $(1,0,0)$ in $S^2$.  One can easily check that the
set of points
$$C=\{(\cos t, \sin t, 0,0)\; | \; t \in \R\}$$ in $S^3$ all map to
$(1,0,0)$ via the Hopf map $h$.  In fact, this set $C$ is the {\em
entire} set of points that map to $(1,0,0)$ via $h$.  In other words,
$C$ is the preimage set $h^{-1}((1,0,0))$.  You may recognize that $C$
is the unit circle in a plane in $\R^4$.  As we shall see, this is
typical: for {\em any} point $P$ in $S^2$, the preimage set $h^{-1}(P)$
is a circle in $S^3$.  We will also refer to the preimage set
$h^{-1}(P)$ as the {\em fiber} of the Hopf map over $P$.

We devote the remainder of this article to study one aspect of the
geometry of the Hopf fibration, namely, the configuration of its fibers
in $S^3$.  Using stereographic projection (to be explained below) we get
a particularly elegant decomposition of 3-space into a union of disjoint
circles and a single straight line.  Because this arrangement is fun to
think about, we cast it first in the form of a puzzle.

\begin{rem}{Linked Circles Puzzle}{circfolpuzz}
Using disjoint circles and a single straight line, can you fill up 3-space
in such a way that each pair of circles is linked, and the line passes
through the interior of each circle?
\end{rem}

\begin{figure}[h]
\begin{center}
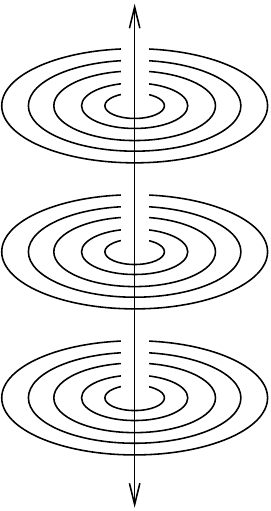
\end{center}
\vspace*{-.2in}\caption{\small One way to fill $\R^3$ with disjoint
circles and a line.  Now try to arrange for every pair of circles to be linked!}
\label{concentric}
\end{figure}

It is the linked-ness of the circles that makes this puzzle
interesting.  If the circles are not required to be linked, there are
easy solutions.  For example, just take stacks of concentric circles
whose centers lie on the given line (see Figure~\ref{concentric}).  We will show
that the Hopf fibers themselves give rise to a solution to this puzzle,
but see if you can think of your own solution first!

We begin with an observation, presented in the form of an Investigation,
on how to find rotations that take a given point $A$ to a given point
$B$.

\begin{figure}[h]
\begin{center}
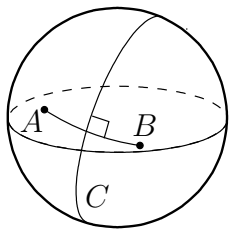
\end{center}
\vspace*{-.2in}\caption{\small The axis of any rotation taking $A$ to $B$ must pass
through the great circle $C$ that bisects $\overline{AB}$.}
\label{bisectgrcirc}
\end{figure}

\begin{figure}[h]
\begin{center}
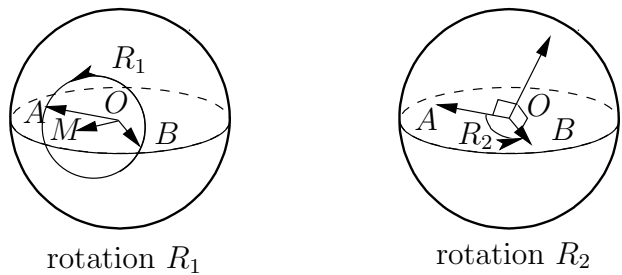
\end{center}
\vspace*{-.2in}\caption{\small Two rotations taking $A$ to $B$.}
\label{specialrots}
\end{figure}

\begin{rem}{Investigation}{findrots}
Given two points $A$ and $B$ on $S^2$, how can we describe the set of
{\em all} possible rotations that move $A$ to $B$?  First, choose an arc
of a great circle joining $A$ to $B$ and call this arc $\overline{AB}$;
note that this choice is not unique.  Convince yourself that if $R$ is a
rotation taking $A$ to $B$, then the axis of $R$ must lie somewhere
along the great circle bisecting $\overline{AB}$ (see
Figure~\ref{bisectgrcirc}).

Along this great circle there are two axes of rotation for which
the angle of rotation is easy to compute.
\begin{enumerate}
\item \label{midrot} When the axis of rotation
passes through the midpoint $M$ of $\overline{AB}$, the angle of
rotation $\theta$ is $\pi$ radians or 180 degrees.  Let us call this
rotation $R_1$  (see the drawing on the left in Figure~\ref{specialrots}).
\item \label{perprot} When the axis of
rotation is perpendicular to the vectors ${\mathbf v}=\vec{OA}$ and
${\mathbf w}=\vec{OB}$, the angle of
rotation is (plus or minus) the angle between ${\mathbf v}$ and
${\mathbf w}$ and is given
by $\cos(\theta)={\mathbf v}\cdot{\mathbf w}$.  We will call this
rotation $R_2$  (see the drawing on the right in Figure~\ref{specialrots}).
\end{enumerate}
\end{rem}

If a point $r$ in $S^3$ is sent by the Hopf map to the point $P$
in $S^2$, then by~(G)
we know that the rotation $R_r$
moves the point $(1,0,0)$ to $P$.  We can use~(I)
to find the axis and angle of rotation for two rotations that map
$(1,0,0)$ to $P$.  Let us call $R_1$ and $R_2$ the rotations described
in~(I)
parts~\ref{midrot} and~\ref{perprot},
respectively. 

Once we have axes and angles of rotation for $R_1$ and $R_2$, we can
use~(D)
to find the quaternions $r_1$ and $r_2$ that
map to $R_1$ and $R_2$ under the map $\varphi$, i.e., $R_1=R_{r_1}$ and
$R_2=R_{r_2}$. 

\begin{rem}{Investigation}{fibdescr}
What are explicit formulas for the quaternions $r_1$ and $r_2$ described
above?  For the point $P=(p_1,p_2,p_3)$ on $S^2$, verify that the quaternions $r_1$
and $r_2$ are given by
\begin{eqnarray*}
r_1 &=&
\frac{1}{\sqrt{2(1+p_1)}}\left((1+p_1)i+p_2j+p_3k\right),\\
r_2 &=&
\sqrt{\frac{1+p_1}{2}}\left(1+\frac{-p_3j}{1+p_1}+\frac{p_2k}{1+p_1}\right).
\end{eqnarray*}
Let us write $e^{it}$ for $\cos t + i\sin t$.
The fiber $h^{-1}(P)$ is given as a parametrically defined circle in
$\R^4$ by either of the following.
\begin{eqnarray*}
h^{-1}(P) &=& \{r_1e^{it}\}_{0\leq t\leq 2\pi}\\
h^{-1}(P) &=& \{r_2e^{it}\}_{0\leq t\leq 2\pi}
\end{eqnarray*}
The point $P=(-1,0,0)$ is a special case, and $h^{-1}((-1,0,0))$ is
given by
$$h^{-1}((-1,0,0))=\{ke^{it}\}_{0\leq t\leq 2\pi}.$$
\end{rem}

\section*{Seeing the Hopf fibration}

Next we demonstrate a method that allows us to {\em see} a little of
what is going on with the Hopf fibration.  Our aim is to show pictures
of fibers.  We do this by means of {\em stereographic projection}.

\begin{figure}[h]
\begin{center}
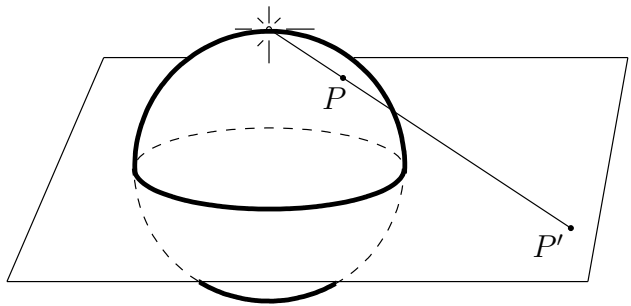
\end{center}
\vspace*{-.2in}\caption{\small Stereographic projection.}
\label{stereoproj}
\end{figure}

We begin by describing the stereographic projection of the 2-sphere to
the $x,y$-plane (see Figure~\ref{stereoproj}).  Imagine a light source placed
at the ``north pole'' $(0,0,1)$.  Stereographic projection sends a point
$P$ on $S^2$ to the intersection of the light ray through $P$ with the
plane.

The alert reader will notice that the point $(0,0,1)$ has no sensible
image under this projection.  Therefore we restrict the stereographic
projection to $S^2\setminus(0,0,1)$.

\begin{rem}{Investigation}{stereoprojection}
Verify that the stereographic projection described above is given by
$$(x,y,z)\mapsto \left(\frac{x}{1-z},\frac{y}{1-z}\right).$$ Write out
the {\em inverse} map $\R^2 \rightarrow S^2\setminus (0,0,1)$.  That is,
given a point $(a,b)$ in the plane, what are the $(x,y,z)$ coordinates
of the point on $S^2$ sent to $(a,b)$ by the stereographic projection?
Show that a circle on $S^2$ that contains $(0,0,1)$ is mapped to a
straight line in the plane.  Prove that a circle on $S^2$ that does not
pass through the point of projection $(0,0,1)$ is mapped by the
stereographic projection to a circle in the plane.  (For a proof of the
circle preservation property using elementary geometry of complex
numbers, see~\cite{ahlfors} Ch.~1 \S~2.4.)
\end{rem}

Like the definition of the sphere, stereographic projection generalizes
to all dimensions, and in particular, it provides a projection map
$S^3\setminus (1,0,0,0) \rightarrow\R^3$ given by
\begin{equation}\label{stereothree}
(w,x,y,z)\mapsto
\left(\frac{x}{1-w},\frac{y}{1-w},\frac{z}{1-w}\right).
\end{equation}
Note that the point $(1,0,0,0)$ on $S^3$ from which we project is an arbitrary choice.

The real power of stereographic projection is this: it allows us to see
all of the 3-sphere (except one point) in familiar 3-space.  This is
remarkable because $S^3$ is a curved object that resides in 4-space.  

The last property in~(K)
above---that stereographic
projection preserves circles---holds in all dimensions
(see~\cite{berger}, Chapter~18).
We know from the previous section that fibers of the Hopf map are
circles in $S^3$.  It follows that stereographic projection sends them
to circles (or a line, if the fiber contains the point $(1,0,0,0)$) in
$\R^3$.  We conclude with two Investigations that show how the
stereographic images of the Hopf fibers solve the linked circles
puzzle~(H).

\begin{figure}[h]
\begin{center}
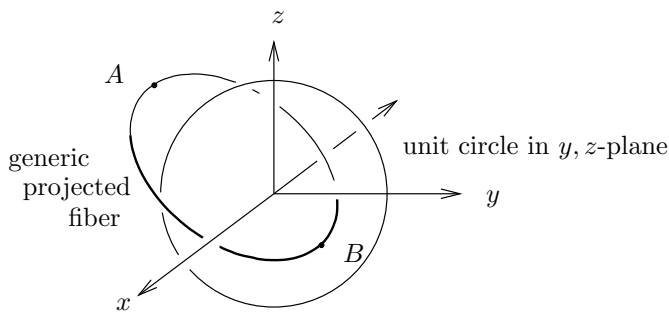
\end{center}
\vspace*{-.2in}\caption{\small A generic projected Hopf fiber. $A$ and $B$ mark the
intersections of the fiber with the $y,z$-plane.}
\label{genericfib}
\end{figure}

\begin{rem}{Investigation}{puzzsoln1}
Let us denote by $s$ the stereographic projection $S^3\setminus
(1,0,0,0)\rightarrow \R^3$ given in~(\ref{stereothree}).  Then $s\of
h^{-1}((1,0,0))$ is the $x$-axis, $s\of h^{-1}((-1,0,0))$ is the unit
circle in the $y,z$-plane, and for any other point $P=(p_1,p_2,p_3)$ on
$S^2$ not equal to $(1,0,0)$ or $(-1,0,0)$, $s\of h^{-1}(P)$ is a circle
in $\R^3$ that intersects the $y,z$-plane in exactly two points $A$ and
$B$, one inside and one outside the unit circle in the $y,z$-plane.
This establishes that $s\of h\inv(P)$ is linked with the unit circle in
the $y,z$-plane.  The points $A$ and $B$ lie on a line through the
origin containing the vector $(0,p_3,-p_2)$.  The plane of the circle
$s\of h\inv(P)$ cannot contain the $x$-axis (if it did, $s\of h\inv(P)$
would intersect $s\of h\inv((1,0,0))$, but fibers are disjoint).  From
these observations we can conclude that the $x$-axis passes through the
interior of the circle $s\of h\inv(P)$.  See Figure~\ref{genericfib}.
\end{rem}

\begin{figure}[h]
\begin{center}
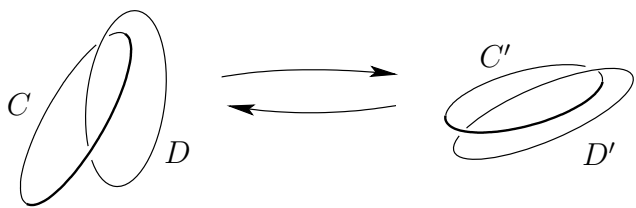
\end{center}
\vspace*{-.2in}\caption{\small If the continuous bijective images $C'$,
$D'$ of circles $C$, $D$ are linked, then $C$ and $D$ must also be
linked.}
\label{pairolinks}
\end{figure}

\begin{rem}{Investigation}{puzzsoln2}
To show that {\em any two} projected fiber circles $C$ and $D$ are
linked, we exhibit a continuous one-to-one map $\psi\colon\R^3
\rightarrow \R^3$ that takes $C$ to the unit circle in the $y,z$-plane,
and takes $D$ to some other projected fiber circle $E$.  Since $E$ is
linked with the unit circle in the $y,z$-plane, $C$ and $D$ must also be
linked.  See Figure~\ref{pairolinks}.  {(Students who have never studied
topology may accept the intuitively reasonable statement that
linked-ness of circles cannot be altered by a continuous bijective map.
Students with experience in topology may enjoy trying to prove this.)}

Here is how to construct the map $\psi$.  Let $P$ be any point on the
circle $C$, and let $r=s^{-1}(P)$.  Define $f\colon\R^4\rightarrow\R^4$
by $f(x)=kr^{-1}x$ (quaternion multiplication).  The map $\psi$ is the
composition $s\of f \of s^{-1}$.
\end{rem}

\begin{figure}[h]
\begin{center}
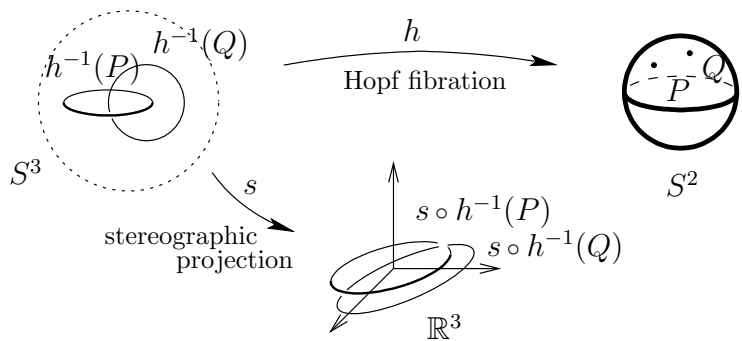
\end{center}
\vspace*{-.2in}\caption{\small Stereographic projections of Hopf fibers.  Any two
projected fibers are linked circles, except $s\of h^{-1}(1,0,0)$ is a line.}
\end{figure}

\subsection*{Conclusion}

We have explained how the Hopf fibration can be understood in terms of
quaternions.  In the process, we showed how the algebra of
rotations in 3-space is built into the workings of the Hopf map.

\begin{figure}[h]
\begin{center}
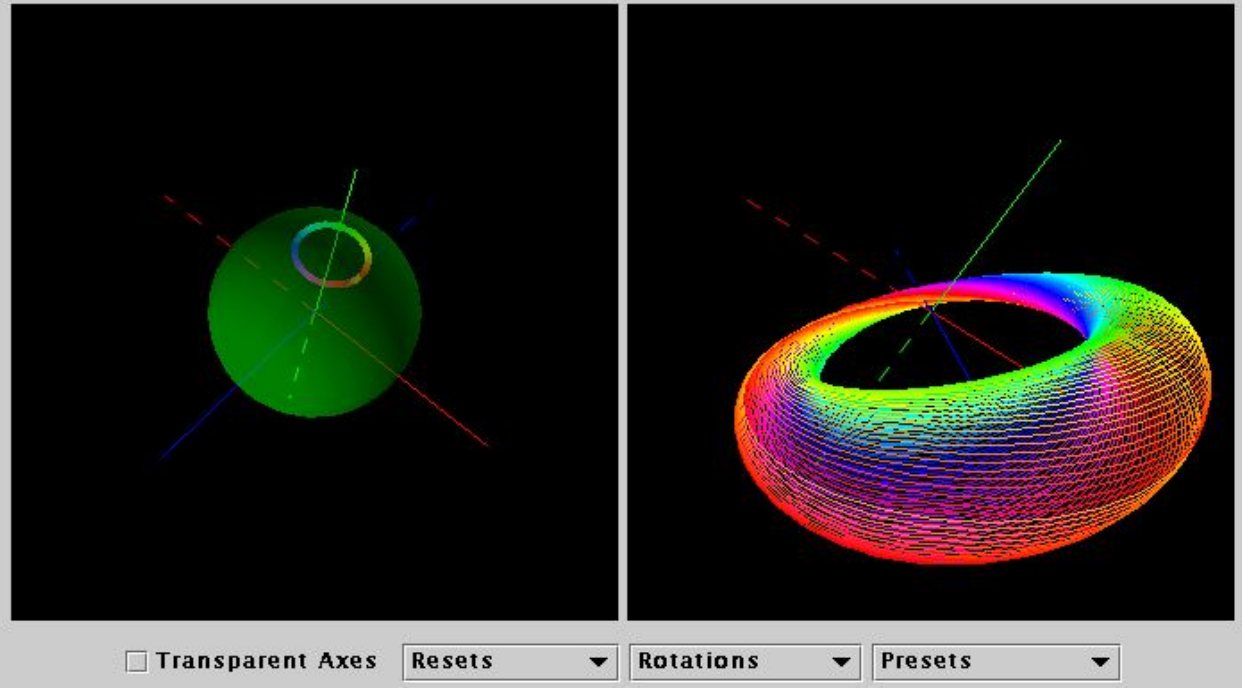
\end{center}
\vspace*{-.2in}\caption{\small Screenshot of Hopf fiber software.}
\label{sh}
\end{figure}

Topics raised in the Investigations suggest many lines of inquiry for
independent student research.  For example, making computer animations
of linked Hopf fibers has been an independent study research project for
two of our undergraduate students.  Figure~\ref{sh} shows an image from
the software written by Nick Hamblet (see Acknowledgment below).  The
left panel shows a set of points lying on a circle in the codomain $S^2$
of the Hopf fibration.  The right panel shows, via stereographic
projection, the fibers corresponding to those points.  An ongoing
project is to build a web tutorial site featuring the animations.  The
reader who finds topics in this article appealing will enjoy a related
article~\cite{zulli}.  For general inspiration, and more on the geometry
of $\R^3$ and rotations, see Hermann Weyl's lovely book {\it
Symmetry}~\cite{weyl}.

\paragraph*{Acknowledgment.}
We are grateful to Lebanon Valley College for summer support for Nick
Hamblet's software development project.  Nick Hamblet is a student at
Lebanon Valley College, class of 2004.  Nick's work continues a similar
project at Wake Forest University with Paul Hemler, Professor of Computer
Science, and Keely Chorn, class of 1999.

\end{document}